\documentclass[12pt]{amsart}

\usepackage{amssymb}

\usepackage[latin1]{inputenc}
\usepackage{mathptmx}
\usepackage[all]{xy}
\theoremstyle{plain}
  \newtheorem{theorem}{Theorem}[section]
  \newtheorem{proposition}[theorem]{Proposition}
  \newtheorem{lemma}[theorem]{Lemma}
  \newtheorem{corollary}[theorem]{Corollary}
\theoremstyle{definition}
  \newtheorem{definition}[theorem]{Definition}
  \newtheorem{example}[theorem]{Example}
  
  \newtheorem{remark}[theorem]{Remark}
  
\newtheorem{conjecture}[theorem]{Conjecture}

\newcommand{\Z}{\mathbb{Z}}

\newcommand{\N}{\mathbb{N}}

\newcommand{\Ac}{\mathcal{A}}
\newcommand{\Bc}{\mathcal{B}}

\newcommand{\Mcc}{\mathcal{M}}

\newcommand{\vb}{{\bf v}}

\newcommand{\cx}{\mathrm{cx}}

\newcommand{\KER}{\mathrm{Ker}}
\newcommand{\revlex}{\mathrm{rlex}}
\newcommand{\lex}{\mathrm{lex}}

\DeclareMathOperator{\depth}{depth}

\DeclareMathOperator{\GL}{GL}

\DeclareMathOperator{\gin}{gin}

\DeclareMathOperator{\ini}{in}

\DeclareMathOperator{\projdim}{proj\,dim}

\DeclareMathOperator{\reg}{reg}

\DeclareMathOperator{\Span}{span}

\DeclareMathOperator{\supp}{supp}

\DeclareMathOperator{\Tor}{Tor}

\DeclareMathOperator{\tensor}{\otimes}

\DeclareMathOperator{\pnt}{\raise 0.5mm \hbox{\large\bf.}}

\DeclareMathOperator{\mlpnt}{\!\!\hbox{\large\bf.}}

\begin{document}
\title{Exterior depth and exterior generic annihilator numbers}

\author{Gesa K\"ampf}
\address{Fachbereich Mathematik und Informatik\\
      Universit\"at Osnabr\"uck\\
      49069 Osnabr\"uck, Germany}
\email{gkaempf@mathematik.uni-osnabrueck.de}

\author{Martina Kubitzke}
\address{Fachbereich Mathematik und Informatik\\
      Philipps-Universit\"at Marburg\\
      35032 Marburg, Germany}
\email{kubitzke@mathematik.uni-marburg.de}

\begin{abstract}
We study the exterior depth of an $E$-module and its exterior generic annihilator numbers. For the exterior depth of a squarefree $E$-module we show how it relates to the symmetric depth of the corresponding $S$-module and classify those simplicial complexes having a particular exterior depth in terms of their exterior algebraic shifting. 
We define exterior annihilator numbers analogously to the annihilator numbers over the polynomial ring introduced by Trung and Conca, Herzog and Hibi. In addition to a combinatorial interpretation of the annihilator numbers we show how they are related to the symmetric Betti numbers and the Cartan-Betti numbers, respectively. We finally conclude with an example which shows that neither the symmetric nor the exterior generic annihilator numbers are minimal among the annihilator numbers with respect to a sequence.
\end{abstract}

\maketitle

\section{Introduction}
Let $S=K[x_1,\ldots,x_n]$ be the symmetric algebra and $E=K\langle e_1,\ldots,e_n\rangle$ denote the exterior algebra over an infinite field $K$.
Let $\Delta$ be a simplicial complex on the vertex set $[n]=\{1,\ldots,n\}$.  Instead of studying the properties of $\Delta$ directly one often studies the \emph{Stanley-Reisner ring} $K[\Delta]=S/I_\Delta$,  or the \emph{exterior face ring} $K\{\Delta\}=E/J_\Delta$  of $\Delta$. Several algebraic and homological invariants of $K[\Delta]$ over $S$ are analogous to invariants of $K\{\Delta\}$ over $E$. In this paper we study relations between the corresponding symmetric and exterior invariants, also in the more general context of squarefree modules.

Yanagawa introduced squarefree modules over the polynomial ring as a generalization of squarefree monomial ideals in \cite[Definition 2.1]{Y}. 
R\"omer defined in \cite[Definition 1.4]{R} the corresponding notion of a squarefree $E$-module and showed that there exists an equivalence of categories between the category of squarefree $S$-modules and the category of squarefree $E$-modules (see Section \ref{exterior depth} for details). 
A typical example of a squarefree $S$-module is the Stanley-Reisner ring $K[\Delta]$, 
of a squarefree $E$-module the exterior face ring $K\{\Delta\}$.
Under the equivalence of categories $K[\Delta]$ corresponds to $K\{\Delta\}$. 

Let $\Mcc$ be the category of finitely generated graded left and right $E$-modules $M$ satisfying $am=(-1)^{\deg a\deg m}ma$ for homogeneous elements $a\in E$, $m\in M$.
For example, if $J\subseteq E$ is a graded ideal, then $E/J$ belongs to $\Mcc$.

Following \cite{AAH} we call a linear form $v\in E_1$  \emph{regular} on $M\in\Mcc$ if the annihilator of $v$ is the smallest possible, i.e., if $0:_M v=v\cdot M$. 
A regular sequence on $M$ is a sequence $v_1,\ldots,v_r$ in $E_1$ such that $v_i$ is regular on $M/( v_1,\ldots,v_{i-1}) M$ for $1\leq i\leq r$ and $M/( v_1,\ldots,v_{r}) M\not=0$. 
The \emph{exterior depth} of $M$ is also introduced in \cite{AAH} as the length of a maximal regular sequence and denoted by $\depth_E M$.
In the second section of this paper we study the relation between the depth of a squarefree $S$-module and the exterior depth of the corresponding squarefree $E$-module. 

 Using some results of Aramova, Avramov and Herzog \cite{AAH} and R\"omer \cite{R} we are able to show in Theorem \ref{ineqDepths} that for a squarefree $S$-module $N$ and its corresponding squarefree $E$-module $N_E$ it holds that
\begin{equation}\label{ineq:intro1}
\depth_E(N_E)\leq \depth_S(N).\end{equation} 
As byproduct we additionally obtain the inequalities
\begin{equation}\label{ineq:intro3}
\cx_E (N_E)\geq \projdim_S(N) \qquad \text{ and } \qquad
\depth_S(N)-\depth_E(N_E)\leq \reg_S (N),\end{equation} 
where $\cx_E (N_E)$ is the complexity of $N_E$ (see Section \ref{exterior depth} for a definition).

It is an interesting question to classify all modules for which equality holds in the second inequality of (\ref{ineq:intro3}). We show in Lemma \ref{CM} that this is the case for $K[\Delta]$ if $\Delta$ is a Cohen-Macaulay simplicial complex or if its Stanley-Reisner ideal has a linear resolution. But one easily finds examples of non-Cohen-Macaulay complexes for which equality holds even though their Stanley-Reisner ideal has no linear resolution (cf. Example \ref{ex:notCM}). 

Starting with a simplicial complex $\Delta$ one can construct its so-called exterior shifting $\Delta^e$ (see Section \ref{exterior depth} for more details). This passage corresponds to the transition of an ideal in the exterior algebra to its generic initial ideal. It allows to compute some invariants defined in terms of an ideal and a generic basis by computing the correponding invariants for the generic initial ideal and the standard basis. This idea goes back a long way, see e.g., \cite{AH}, \cite{AHH3}, \cite{BK}, \cite{HT}, \cite{Kalai}.

The exterior depths of $K\{\Delta\}$ and $K\{\Delta^e\}$ coincide as was shown in \cite[Proposition 2.3]{HT}. 
In Theorem \ref{characterization} we describe the structure of the exterior shifted complex of $\Delta$ in terms of the exterior depth of $K\{\Delta\}$, i.e., we show that 
$\depth_E (K\{\Delta\})=r$ if and only if $\Delta^{e}=2^{[r]}\ast\Gamma$. 
Here, $\Gamma$ is a non-acyclic simplicial complex, $2^{[r]}$ denotes the $(r-1)$-simplex. 

Finally, we construct for every triple $(s,t,r)$ of natural numbers with $r\geq s-t\geq 0$ a simplicial complex $\Delta$ such that $\depth_E(K\{\Delta\})=s$, $\depth_S(K[\Delta])=t$ and $\reg_S(K[\Delta])=r$.

In the third section our discussion is concentrated on the \emph{exterior annihilator numbers} $\alpha_{i,j}(v_1,\ldots,v_n;M)$ of an $E$-module $M$ with respect to a sequence $v_1,\ldots,v_n$ of linear forms. These numbers are the exterior analogue of the symmetric annihilator numbers, first defined by Trung in \cite{Trung}
and subsequently studied by Conca, Herzog and Hibi in \cite{CHH}. The results in this section are very much in the spirit of their results. 
We show that there is a non-empty Zariski-open set on which the annihilator numbers are constant.
This gives rise to the definition of the \emph{exterior generic annihilator numbers} $\alpha_{i,j}(E/J)$ for a graded ideal $J$ in $E$. 
These numbers are a refinement of the exterior depth. 
For the exterior face ring of a simplicial complex  they have a nice combinatorial interpretation, namely (see Corollary \ref{description annihilators of complex}) 
$$\alpha_{i,j}(K\{\Delta\})=|\{F\in\Delta^e~:~ |F|=j,\  [i]\cap F=\emptyset,\  F\cup\{i\}\not\in\Delta^e\}|.$$

Using this combinatorial description we can express in Proposition \ref{relation with polynomial ring} the symmetric Betti numbers of $K[\Delta^e]$ as a linear combination of certain exterior generic annihilator numbers, more precisely
$$\beta^S_{i,i+j}(K[\Delta^e])=\sum_{l=1}^n \binom{n-l-j}{i-1}\alpha_{l,j}(K\{\Delta\}).$$

Finally, we show in Theorem \ref{CartanBetti} that the so-called Cartan-Betti numbers (see \cite[Definition 2.2]{NRV}) can be bounded from above by a positive linear combination of the exterior generic annihilator numbers. In particular, this gives rise to an upper bound for the exterior Betti numbers, 
$$\beta^E_{i,i+j}(E/J)\leq\sum_{k=1}^n \binom{n+i-k-1}{i-1}\alpha_{k,j}(E/J) \quad i\geq 1, j\geq 0,$$
which has been also obtained in another way in \cite[Theorem 2.4(i)]{NRV}.
Furthermore, this bound is tight if and only if $J$ is componentwise linear. 

In the fourth section we discuss the issue if the generic annihilator numbers are the minimal ones among all annihilator numbers with respect to a sequence. By minimal we mean that for any basis $\{v_1,\ldots,v_n\}$ of 
$E_1$, $1\leq i\leq n$ and $j\in \mathbb{Z}$ it should hold that
$$\alpha_{i,j}(E/J)\leq \alpha_{i,j}(v_1,\ldots,v_n;E/J).$$
We give an example, where this is not the case (see Example \ref{not minimal}). After a slight modification the same example shows that also in the symmetric case this minimality is not given (see Example \ref{not minimal 2}).

\section{Exterior depth}\label{exterior depth}

Yanagawa introduced squarefree modules over the polynomial ring as a generalization of squarefree monomial ideals in \cite{Y}.
We fix some notations.
For $a=(a_1,\ldots,a_n)\in\N^n$ we say that $a$ is \emph{squarefree} if $0\leq a_i \leq 1$ for $i=1,\ldots,n$. We set $|a|=a_1+\ldots+a_n$ and $\supp(a)=\{i~ :~a_i\neq 0\}$.

A finitely generated $\mathbb{N}^n$-graded $S$-module $N=\oplus_{a\in\mathbb{N}^n}N_a$ is called  \emph{squarefree} if the multiplication map $N_a\rightarrow N_{a+\epsilon_i}
,\hspace{5pt}y\mapsto x_i y$ is bijective for any $a\in\mathbb{N}^n$ and for all $i\in\supp(a)$, where $\epsilon_i\in\N^n$ is the vector with $1$ at the $i$-th position and zero otherwise.

R\"omer defined in \cite[Definition 1.4]{R} a finitely generated $\mathbb{N}^n$-graded $E$-module $M=\oplus_{a\in\mathbb{N}^n}M_a$ to be \emph{squarefree} if it has only squarefree (non-zero) components. 

Aramova, Avramov and Herzog and R\"omer construct in \cite{AAH} and \cite{R} a squarefree $E$-module $N_E$ and its minimal free resolution starting from a squarefree $S$-module $N$ and its minimal free resolution.

\begin{theorem}\cite[Theorem 1.2]{R}
The assignment $N\mapsto N_E$ induces an equivalence between the categories of squarefree $S$-modules and squarefree $E$-modules (where the morphisms are the $\N^n$-graded homomorphisms). 
\end{theorem}

For $F=\{i_1,\ldots,i_r\}\subseteq[n]$ we set $x_F=x_{i_1}\cdot\ldots\cdot x_{i_r}$ and $e_F=e_{i_1}\wedge\ldots\wedge e_{i_r}$. We usually assume that $1\leq i_1<\ldots<i_r\leq n$. The elements $e_F$ are called \emph{monomials} in $E$.

\begin{example}
 Let $\Delta$ be a simplicial complex on the vertex set $[n]$. Recall that a simplicial complex $\Delta$ on vertex set $[n]$ is a collection of subsets of $[n]$ such that whenever $F\in \Delta$ and $G\subseteq F$ it holds that $G\in \Delta$. The elements of $\Delta$ are called \emph{faces}. Throughout this paper we always assume that $i\in \Delta$ for all $i\in [n]$.  

Let $K[\Delta]=S/I_\Delta$ be the \emph{Stanley-Reisner ring} of $\Delta$, where $I_\Delta$ is the Stanley-Reisner ideal $I_\Delta=(x_F:F\not\in\Delta)$ of $\Delta$, and let $K\{\Delta\}=E/J_\Delta$ be the \emph{exterior face ring}, where $J_\Delta$ is the exterior face ideal $J_\Delta=(e_F:F\not\in\Delta)$ of $\Delta$.

The Stanley-Reisner ring $K[\Delta]$ is a typical example of a squarefree $S$-module (see \cite{Y}), the exterior face ring $K\{\Delta\}$ of a squarefree $E$-module (see \cite{R}).

These two correspond to each other under the equivalence of categories, i.e. $K[\Delta]_E\cong K\{\Delta\}$ (see \cite[Theorem 1.3]{AAH}).
\end{example}

Our aim is to show that the symmetric depth of a squarefree $S$-module $N$ is always greater or equal than the exterior depth of the corresponding squarefree $E$-module $N_E$.  Recall from the introduction that the depth of an $E$-module $M\in\Mcc$ is the length of a maximal $M$-regular sequence.
It is strongly related with the \emph{complexity} $\cx_E (M)$ of the module which is defined as 
$$\cx_E (M)=\inf\{c\in \mathbb{Z}~:~\beta_i^{E}(M)\leq \alpha i^{c-1}\text{ for some }\alpha\in\mathbb{R}\text{ and for all }i\geq 1\}.$$
In other words, the complexity measures the polynomial growth of the exterior Betti numbers of $M$ and is therefore a measure for the size of a minimal free resolution of $M$ by free $E$-modules.

The proof is based on comparisons of several invariants of modules over $S$ or $E$. We collect the formulas used for these comparisons in the next theorem.

\begin{theorem}\label{thm:formulas}
\begin{enumerate}
\item[(i)] \cite[Corollary 1.3]{R} 
  Let $N$ be a squarefree $S$-module and let $N_E$ be the associated squarefree $E$-module. Let $\beta^{E}_{i,j}$ and $\beta^{S}_{i,j}$ denote the graded Betti numbers over $E$ and $S$, respectively. Then
     $$\beta_{i,i+j}^{E}(N_E)=\sum_{k=0}^{i}\binom{i+j-1}{ j+k-1}\beta_{k,k+j}^{S}(N).$$
\item[(ii)] $($e.g., \cite[Theorem 1.3.3]{BH-book}$)$ 
  Let $N$ be a finitely generated $S$-module. Then
     $$ \projdim_S(N)+\depth_S(N)=n.$$
\item[(iii)] \cite[Theorem 3.2]{AAH} 
  Let $M\in\Mcc$. Then
     $$\depth_E(M)+\cx_E (M)=n.$$
\end{enumerate}
\end{theorem}

We have now provided all notions and facts we need to prove the desired inequality.

\begin{theorem}\label{ineqDepths}
Let $N$ be a finitely generated $\mathbb{N}^n$-graded squarefree $S$-module. Then
$$0\leq \depth_S(N)-\depth_E(N_E)=\cx_E (N_E)-\projdim_S(N)\leq \reg_S (N).$$
\end{theorem}
\begin{proof}
A comparison of Theorem \ref{thm:formulas} (ii) and (iii) shows that the differences between the depths and between the complexity and the projective dimension are equal. We show that the inequalities are satisfied by the latter difference. 

From Theorem \ref{thm:formulas}(i) it follows that
$$\beta_i^{E}(N_E)=\sum_{j\geq 0}\beta_{i,i+j}^{E}(N_E)=\sum_{j\geq 0}\sum_{k=0}^{i}\binom{i+j-1}{ j+k-1}\beta_{k,k+j}^{S}(N).$$
Let $m_j^{(i)}=\max\{k+j~:~\beta_{k,k+j}^{S}(N)\neq 0,\; 0\leq k\leq i\}$. From the above formula for the Betti numbers we conclude that $\beta_{i,i+j}^{E}(N_E)$ is a polynomial in $i$ of degree $m_j^{(i)}-1$. Therefore, it follows that $\beta_i^{E}(N_E)$ is a polynomial in $i$ of degree $m^{(i)}-1$, where $m^{(i)}=\max\{m_j^{(i)}~:~j\geq 0, \beta^E_{i,i+j}\not=0\}$ (here we can write max instead of sup because the modules in the minimal free $E$-resolution of a finitely generated module are finitely generated). This yields for the complexity
\begin{eqnarray*}
\cx_E (N_E)&=&\sup\{m^{(i)}~:~i\geq 0\}\\
&=&\sup\{k+j~:~\beta_{k,k+j}^{S}(N)\neq 0,\; k\geq 0,\; j\geq 0\}\\
&=&\max\{l~:~\beta_{i,l}^{S}\neq 0,\; i\geq 0,\;l\geq 0\},
\end{eqnarray*}
where the last equality holds because $N$ has a finite $S$-resolution. Let $p=\projdim_S (N)$. Then there exists $0\leq k\leq\reg_S(N)$ such that $\beta_{p,p+k}^{S}(N)\neq 0$. This implies $\cx_E(N_E)\geq p+k\geq p=\projdim_S (N)$. This finally shows the claim.
\end{proof}

An interesting question to ask is if there are classes of squarefree modules for which equality holds in the second inequality in Theorem \ref{ineqDepths}. In the special case of Stanley-Reisner rings of simplicial complexes we can identify at least two such classes.

\begin{lemma} \label{CM}
Let $\Delta$ be a simplicial complex. If $J_\Delta$ has a linear resolution or if $\Delta$ is Cohen-Macaulay, then
$$\depth_S(K[\Delta])-\depth_E(K\{\Delta\})=\cx_E (K\{\Delta\})-\projdim_S(K[\Delta])= \reg_S (K[\Delta]).$$
\end{lemma}

\begin{proof}
Examining the proof of Theorem \ref{ineqDepths} we see that we have an equality if and only if $\beta^ S_{p,p+r}(K[\Delta])\not=0$ where $p=\projdim_S(K[\Delta])$ and $r=\reg_S(K[\Delta])$. This is the most right lower corner in the Betti diagram. It is obviously the case if $I_\Delta$ has a linear resolution.

If $\Delta$ is Cohen-Macaulay, $K[\Delta]$ has depth $d$ over $S$ where $\dim \Delta=d-1$ and the face ring $K\{\Delta\}$ has a $d$-linear injective resolution over $E$ (see e.g., \cite[Corollary 7.6]{AH} or \cite[Example 5.1]{KR}). Following \cite[Theorem 5.3]{KR} it holds that 
$$ d=\depth_E (K\{\Delta\})+\reg_E (K\{\Delta\}).$$
The regularity over $E$ is the same as the regularity over $S$ (this follows from the relation between the Betti numbers, Theorem \ref{thm:formulas}(i)), hence
$$\reg_S(K[\Delta])=d-\depth_E(K\{\Delta\})=\depth_S (K[\Delta])-\depth_E(K\{\Delta\}).$$
As an alternative proof of this case it is possible to show that if $\Delta$ is Cohen-Macaulay then $K[\Delta]$ has only one extremal Betti number. 
\end{proof}

The following example shows that in general the converse of Lemma \ref{CM} is not true.
\begin{example}\label{ex:notCM}
Let $\Delta$ be the simplicial complex on the vertex set $[4]$ consisting of two triangles -- one of them filled, the other one missing -- glued together along one edge. The face ideal is $J_\Delta=(e_3e_4,e_1e_2e_4)$. It is not generated in one degree and thus does not have a linear resolution. On the other hand $\Delta$ is not pure so it cannot be Cohen-Macaulay. But $I_\Delta$ is squarefree stable whence we can compute the depths and the regularity from known formulas (A monomial ideal $I$ in $S$ is called \emph{squarefree stable}, if for all squarefree monomials $x_A\in S$ and all $i>\min(A)$ with $i\not\in A$ one has $x_ix_{A\setminus\min(A)}\in I$, see Example \ref{ex:s t r independent} for the formulas)). Then
$$\depth_S (K[\Delta]) = \min\{\min(u)+\deg(u): u\in G(I_\Delta)\}-2=2,$$
$$\reg_S(K[\Delta])=\max\{\deg(u):u\in G(I_\Delta)\}-1=2,$$
$$\depth_E (K\{\Delta\})=\min\{\min(u): u\in G(J_\Delta)\}-1=0.$$
\end{example}

In the remainder of the section we continue to discuss the case of simplicial complexes. 
Starting with a simplicial complex $\Delta$ one can construct its so-called exterior shifting (see \cite{Kalai} for more details). 
A simplicial complex $\Delta$ on a ground set $[n]$ is \emph{shifted} if for every face $F\in \Delta$, $i\in F$ and $j<i$ we have that $(F\setminus\{i\})\cup\{j\}\in\Delta$. 
(Algebraists would probably prefer $j>i$ here.)
Let now $\Delta$ be a simplicial complex on vertex set $[n]$. Let $K'/K$ be a field extension containing the algebraically independent elements $a_{ij}$, $1\leq i,j\leq n$ over $K$, and let $f_1,\ldots,f_n$ with $f_i=\sum_{j=1}^na_{ij}e_j$ for $1\leq i \leq n$ be a generic basis of $E'_1$ where $E'=K'\tensor_K E$. For $A=\{i_1<\ldots<i_r\}\subseteq [n]$ we set $f_A=f_{i_1}\wedge\ldots\wedge f_{i_r}\in E$ and by $\bar{f_A}$ we denote the image of $f_A$ in $K'\{\Delta\}$. Let further $<_{\lex}$ denote the lexicographic order on subsets of $\N$ of the same size, i.e., $A<_{\lex}B$ if and only if $\min((A\setminus B)\cup(B\setminus A))\in A$. In order to define the exterior shifting of $\Delta$ we define the shifting of a family of sets of equal cardinality. We set
$$\Delta^{e}_i=\{A\in \binom{[n]}{ i}~:~\bar{f_A}\notin \Span\{\bar{f_{A'}}~|~A'<_{\lex} A\}\}$$
for $0\leq i\leq \dim\Delta+1$. Finally, the exterior shifting of $\Delta$ is the following simplicial complex $\Delta^{e}=\bigcup_{i=0}^{\dim\Delta+1}\Delta^{e}_i$. Observe that $\Delta^{e}$ may depend on $K$, but not on $K'$. 
One can show that $\Delta^{e}$ is indeed a shifted simplicial complex having the same $f$-vector as $\Delta$. Using the definition of the generic initial ideal of an ideal $J\subseteq E$ one can easily show that $J_{\Delta^{e}}=\gin_{<_{\revlex}}(J_{\Delta})$.
Here $<_{\revlex}$ denotes the reverse lexicographic order with respect to $e_1<\ldots <e_n$. We therefore can compute the exterior shifting of a simplicial complex  by computing the generic initial ideal of the exterior face ideal $J_{\Delta}$.

The following remark shows that the exterior depth is unchanged by exterior shifting. 
\begin{remark} \cite[Proposition 2.3]{HT} \label{Rem:extDepth}
Let $J\subseteq E$ be a graded ideal. Then 
$$\depth_E (E/J)=\depth_E(E/\gin_{<_{\revlex}}(J)),$$
where $\gin_{<_{\revlex}}(J)$ denotes the generic inital ideal of $J$ with respect to $<_{\revlex}$.
\end{remark}

A monomial ideal $J$ in $E$ is called \emph{stable}, if for all monomials $e_A\in J$ and all $i>\min(A)$ with $i\not\in A$ one has $e_ie_{A\setminus\{\min(A)\}}\in J$. 
It is called \emph{strongly stable}, if for all monomials $e_A\in J$ and all $i>j$ with $j\in A,\ i\not\in A$ one has $e_ie_{A\setminus\{j\}}\in J$. 
The face ideal of a simplicial complex is a strongly stable ideal if and only if the complex is shifted. Therefore the generic inital ideal is strongly stable (see also \cite[Proposition 1.7]{AHH} for a direct proof). 
It is well-known that for stable ideals an initial segment of $e_1,\ldots,e_n$ up to the depth is a regular sequence. However, for the convenience of the reader we include a proof of it. 

\begin{lemma}\label{variables regular sequence}
Let $J\subseteq E$ be a stable ideal of depth $t$. Then $e_1,\ldots,e_t$ is a regular sequence on $E/J$. 
\end{lemma}
\begin{proof}
If $t=0$ there is nothing to prove, thus we assume $t>0$. 
We show that $e_1$ is regular on $\gin_{<_{\revlex}}(J)=J$. Then the claim follows by induction on $t$ (note that since $J$ is stable also $J+(e_1)/(e_1)$ is stable in $E/(e_1)\cong K\langle e_2,\ldots,e_{n}\rangle$).

Being a regular sequence is an open condition since it is equivalent to the vanishing of the first Cartan homology with respect to the sequence (see \cite[Remark 3.4]{AAH} for the result and e.g., Section \ref{sect annihilators} for a definition of Cartan homology). Hence there exists an $E/J$-regular sequence $v_1,\ldots,v_t$ of linear forms such that there is a generic automorphism $g$ mapping $e_{i}$ to $v_i$ and $\gin(J)=\ini(g(J))$. 
Since the considered monomial order is the revlex order, we have $\ini(g(J)+e_1)=\ini(g(J))+(e_1)$ and 
$\ini(g(J):e_1)=\ini(g(J)):e_1$ by \cite[Proposition 5.1]{AH} (using the reversed order on $[n]$).
Then the Hilbert functions
\begin{eqnarray*}
&\empty& H(-,E/\gin(J+v_1))=H(-,E/(J+(v_1)))=H(-,E/g(J+(v_1)))\\
&=&H(-,E/(g(J)+e_1))=H(-,E/\ini(g(J)+(e_1)))=H(-,E/(\gin(J)+(e_1)))\end{eqnarray*}
are equal. Analogously one sees $H(-,E/\gin(J:v_1))=H(-,E/\gin(J):e_1)$.
On the other hand $J+(v_1)=J:v_1$ as $v_1$ is $E/J$-regular, hence $H(-,E/\gin(J+(v_1)))=H(-,E/\gin(J:v_1))$. Thus the Hilbert functions of $\gin(J)+(e_1)$ and $\gin(J):e_1$ coincide as well which already implies $\gin(J)+(e_1)=\gin(J):e_1$ because $\gin(J)+(e_1)\subseteq\gin(J):e_1$ is clear. This means that $e_1$ is $E/\gin(J)$-regular. 
\end{proof}

The next theorem characterizes simplicial complexes with a certain exterior depth in terms of the exterior shifting.
For a linear form $v\in E_1$ and $M\in\Mcc$ there is the complex
$$(M,v): \qquad \ldots\longrightarrow M_{j-1}\stackrel{\cdot v}{\longrightarrow}M_j\stackrel{\cdot v}{\longrightarrow}M_{j+1}\longrightarrow \ldots,$$
since $v^2=0$. Then $v$ is $M$-regular if and only if all homology modules $H^j(M,v)$ of the complex vanish for all $j\in\Z$.

If $v\in E_1$ is a generic element then following \cite[Lemma 3.3]{AH} the homology of the complex $(K\{\Delta\},v)$ is isomorphic to the reduced simplicial homology of $\Delta$ (see e.g., \cite[Section 5.3]{BH-book} for a definition of reduced simplicial homology). Therefore, we have $\depth_E (K\{\Delta\})=0$ if and only if there exists $0\leq i\leq \dim\Delta$ such that $\widetilde{H}_i(\Delta;K)\neq 0$. A simplicial complex $\Delta$ with this property is called \emph{non-acyclic}.

\begin{theorem} \label{characterization}
Let $\Delta$ be a simplicial complex on vertex set $[n]$. Then
$\depth_E (K\{\Delta\})=r$ if and only if $\Delta^{e}=2^{[r]}\ast\Gamma$, 
where $\Gamma$ is a non-acyclic simplicial complex, $2^{[r]}$ the $(r-1)$-simplex and $\dim \Gamma=\dim\Delta-r$.
\end{theorem}

\begin{proof}
We first assume that $\depth_E(K\{\Delta\})=r$.
In order to prove the statement we need to show that for $F\in \Delta^{e}$ it holds that $F\cup[t]$ is a face of $\Delta^{e}$ for any $t\leq r$.
Let $F\in \Delta^{e}$. Without loss of generality we may assume that $F\cap [t]=\emptyset$.
We can take $e_1,\ldots, e_t$ as a regular sequence for the exterior face ring $K\{\Delta^{e}\}=E/\gin_{<_{\revlex}}(J_{\Delta})$ by Lemma \ref{variables regular sequence}. 
Suppose $F\in\Delta^{e}$ such that $F\cap [t]=\emptyset$ and $F\cup[t]\notin \Delta^{e}$. 
Thus, $e_t\wedge\ldots\wedge e_1\wedge e_F=0$ in $K\{\Delta^{e}\}$. By the definition of a regular sequence it follows that 
$$e_{t-1}\wedge\ldots\wedge e_1\wedge e_F\in (e_t)\subseteq K\{\Delta^{e}\}.$$
Therefore, $e_{t-1}\wedge \ldots\wedge e_1\wedge e_F=0\in E/(\gin_{<_{\revlex}}(J_{\Delta})+( e_t))$. Inductively, we get $e_F=0\in E/(\gin_{<_{\revlex}}(J_{\Delta})+( e_t,\ldots,e_1))$. Since $F\cap [t]=\emptyset$ by assumption, we have $e_F\in \gin_{<_{\revlex}}(J_{\Delta})$, i.e., $F\notin \Delta^{e}$. This is a contradiction.\\
This shows $\Delta^{e}=2^{[r]}\ast\Gamma$ for some $(\dim\Delta-r)$-dimensional simplicial complex $\Gamma$ with $J_\Gamma=\gin_{<_{\revlex}}(J_\Delta)+(e_1,\ldots,e_r)$. By definition of the exterior depth and Lemma \ref{variables regular sequence} it holds that $\depth_E(K\{\Gamma\})=\depth_E (K\{\Delta^{e}\})-r=0$. We therefore conclude that $\widetilde{H}_i(\Gamma;K)\neq 0$ for some $0\leq i\leq \dim(\Gamma)$, i.e., $\Gamma$ is non-acyclic.\\
In order to prove the sufficiency part recall that Remark \ref{Rem:extDepth} implies that $\depth_E (K\{\Delta\})=\depth_E(K\{\Delta^{e}\})$. Thus, we only need to show that $\depth_E(K\{\Delta^{e}\})=r$. By assumption, we have that $\Delta^{e}=2^{[r]}\ast \Gamma$, where $\Gamma$ is a non-acyclic simplicial complex. Then the sequence $e_1,\ldots,e_r$ is regular on $K\{\Delta^{e}\}$ by Lemma \ref{variables regular sequence}. This implies $\depth_E (K\{\Delta^{e}\})\geq r$. It further holds that $K\{\Delta^{e}\}/( e_1,\ldots,e_r) \cong K\{\Gamma\}$. Since $\Gamma$ is non-acyclic we know from the remarks preceding this theorem that $\depth_E (K\{\Gamma\})=0$, i.e., there does not exist any regular element on $K\{\Gamma\}$. Using that each regular sequence on $K\{\Delta^{e}\}$ can be extended to a maximal one, we therefore deduce that $e_1,\ldots, e_r$ is already maximal and thus it follows that $\depth_E (K\{\Delta^{e}\})=r$. 
\end{proof}

\begin{remark}
The above theorem can be used to deduce two of the previous results in the special case of simplicial complexes.
\noindent
\begin{enumerate}
\item[(i)] Note that from the above proof the inequality 
$$\depth_E (K\{\Delta\})\leq \depth_S (K[\Delta])$$
between the exterior and the symmetric depth can be deduced. 
If $\depth_E (K\{\Delta\})=r$, then Theorem \ref{characterization} implies $\Delta^{e}=2^{[r]}\ast\Gamma$ and therefore $x_1,\ldots,x_r$ is a regular sequence for $K[\Delta^{e}]$, i.e., $\depth_S(K[\Delta^{e}])\geq r$. Since $\depth_S(K[\Delta])=\depth_S(K[\Delta^{e}])$ (this follows from \cite[Theorem 9.7]{AH} which states that the extremal Betti numbers of $K[\Delta]$ and $K[\Delta^{e}]$ coincide whence the projective dimensions and thus the depths by Auslander-Buchsbaum (cf. Theorem \ref{thm:formulas}(ii)) of both are equal) we obtain the required inequality.\\
\item[(ii)] Using the characterization of Theorem \ref{characterization} we can give a second proof of Lemma \ref{CM}. If $\Delta$ is non-acyclic, i.e., $\widetilde{H}_i(\Delta;K)\neq 0$ for some $0\leq i\leq \dim\Delta$, we have $\depth_E (K\{\Delta\})=0$. Since $\Delta$ is Cohen-Macaulay it follows from Reisner's criterion that $\widetilde{H}_i(\Delta;K)=0$ for $i<\dim\Delta$. Thus, $\widetilde{H}_{\dim\Delta}(\Delta;K)\neq 0$. From \cite[Proposition 2.6]{KW} we conclude that $\reg_S (K[\Delta])=\dim\Delta+1$. Combining these two facts and using that $\Delta$ is Cohen-Macaulay, i.e., $\depth_S(K[\Delta])=\dim_S (K[\Delta])$, we obtain
\begin{eqnarray*}
\depth_S(K[\Delta])-\depth_E(K\{\Delta\})&=&\depth_S (K[\Delta])\\
&=&\dim_S K[\Delta]=\dim\Delta+1=\reg_S(K[\Delta]).
\end{eqnarray*}
Let us now assume that $\Delta$ is an acyclic simplicial complex and let $\depth_E(K\{\Delta\})=r$. We may assume that $\Delta=\Delta^{e}$ since $\Delta$ and $\Delta^{e}$ have the same symmetric and exterior depth, respectively, and the same regularity (see \cite[Corollary 1.3]{R}). From Theorem \ref{characterization} we know that $\Delta^{e}=2^{[r]}\ast \Gamma$, where $\Gamma$ is a non-acyclic simplicial complex.
In particular, since $\Delta^{e}$ is Cohen-Macaulay, so is $\Gamma$. This implies $$\depth_S(K[\Gamma])=\dim_S(K[\Gamma])=\dim_S(K[\Delta])-r=\depth_S(K[\Delta])-r.$$
Using that $\depth_E(K\{\Gamma\})=0=\depth_E(K\{\Delta\})-r$ we deduce
$$\depth_S(K[\Delta])-\depth_E (K\{\Delta\})=\depth_S(K[\Gamma])-\depth_E(K\{\Gamma\}).$$
Since $\Gamma$ is a non-acyclic simplicial complex we know from the first part of our considerations that
$\depth_S(K[\Gamma])-\depth_E(K\{\Gamma\})=\reg_S(K[\Gamma])$.
We further have that $K[\Gamma]\cong K[\Delta]/( x_1,\ldots,x_r)$ and $x_1,\ldots,x_r$ is a regular sequence on $K[\Delta]$. Since reducing modulo a regular sequence leaves the regularity unchanged (see e.g., \cite[Proposition 20.20]{Ei} it holds that $\reg_S(K[\Delta^{e}])=\reg_S(K[\Gamma])$. This finally shows the claim.
\end{enumerate}
\end{remark}

A natural question to ask is which triples of numbers $(t,s,r)$ with $r\geq s-t\geq 0$ can occur such that there exists a simplicial complex $\Delta$ with the property that $\depth_E(K\{\Delta\})=t$, $\depth_S (K[\Delta])=s$ and $\reg_S(K[\Delta])=r$. We can answer this issue by showing that all triples of numbers are possible. 

For a mononomial ideal $I$ in the polynomial ring or in the exterior algebra we denote by $G(I)$ the unique minimal monomial generating set of $I$. 
For a monomial $u\in S$ let $\min(u)=\min\{i~:~x_i\mbox{ divides } u\}$ and analogously $\min(u)=\min\{i~:~e_i\mbox{ divides } u\}$ for a monomial $u\in E$.
Recall that a monomial ideal $I\subseteq S$ is called \emph{squarefree stable} if for all squarefree monomials $u\in I$ and all $i>\min(u)$ such that $x_i$ does not divide $u$ one has $x_i\left(\frac{u}{x_{\min(u)}}\right)\in I$.\\

\begin{example}\label{ex:s t r independent}
Let $s,t,r$ be natural numbers with $r\geq s-t\geq 0$. We construct a simplicial complex $\Delta$ with the property that $\depth_E(K\{\Delta\})=t$, $\depth_S (K[\Delta])=s$ and $\reg_S(K[\Delta])=r$ as follows (note that if $r=0$ there must be $i\in[n]$, $\{i\}\not\in\Delta$, violating our assumption made on simplicial complexes in the rest of the paper).

Let $n=t+r+3$ and $\Delta$ be the simplicial complex on the ground set $[n]$ with minimal non-faces \begin{itemize}
\item  $\{n,n-1,\ldots,n-(s-t)+1,i\}$ with $i=n-(s-t),\ldots,t+1$,
\item  $\{n-r-1,n-r,\ldots,n-2,n-1\}$,
\item  $\{n-r-1,n-r,\ldots,\hat{j},\ldots,n-1,n\}$ with $j=n-(s-t)+1,\ldots,n-1$,
\end{itemize}
where $\hat{j}$ means that $j$ is omitted.
By construction $J_\Delta$ is a stable and $I_\Delta$ a squarefree stable ideal. There are the following formulas to compute the depth and the regularity of such ideals (taking into account the reverse order on $[n]$).
\begin{itemize}
 \item[-]
  $\depth_S (K[\Delta]) = \min\{\min(u)+\deg(u): u\in G(I_\Delta)\}-2$ (\cite[Corollary 2.4]{AHH2});
\item[-] 
  $\reg_S(K[\Delta])=\max\{\deg(u):u\in G(I_\Delta)\}-1$ (\cite[Corollary 2.6]{AHH2});
\item[-] 
  $\depth_E (K\{\Delta\})=\min\{\min(u): u\in G(J_\Delta)\}-1$ (\cite[Proposition 3.4]{KR}).
\end{itemize}

We compute these invariants:

\begin{eqnarray*}
 \depth_S (K[\Delta]) & = & \min\{t+1+(s-t+1),n-r-1+(r+1)\}-2\\
                      & = & \min\{s+2,n\}-2\\
                      & = & s
\end{eqnarray*}
because $n=t+r+3\geq s+3$.

$$ \reg_S (K[\Delta])  =  \max\{s-t+1,r+1\}-1 =  r $$
because $r\geq s-t$.

$$ \depth_E (K\{\Delta\}) = \min\{t+1,n-r-1\}-1 = t$$
because $n-r-1=t+2$.
\end{example}

\section{Exterior generic annihilator numbers}\label{sect annihilators}
In this section we introduce the so-called exterior generic annihilator numbers of an 
$E$-module $M$ which are the exterior analogue of the symmetric generic annihilator numbers of an $S$-module $N$ 
introduced first by Trung in \cite{Trung} and studied later by Conca, Herzog and Hibi (see \cite{CHH}). In what follows we derive some facts which relate those numbers to the symmetric Betti numbers of a module.
In order to do so, we need to lay some background and introduce some notation.\\
We first recall the construction of the Cartan complex from \cite[Section 2]{AHH}. It is a very useful complex over $E$ 
which plays a similar role as the Koszul complex for the polynomial ring.
For a sequence
$\vb=v_1,\ldots,v_m$ with $v_i\in E_1$  let
$C\mlpnt(\vb;E)=C\mlpnt(v_1,\ldots,v_m;E)$ be the free divided power
algebra $E\langle x_1,\ldots,x_m\rangle$. It is generated by the
divided powers $x_i^{(j)}$ for $1\leq i\leq m$ and $j\geq 0$ which
satisfy the relations $x_i^{(j)}x_i^{(k)}=((j+k)!/(j!k!))x_i^{(j+k)}$.
Thus $C_i(\vb;E)$ is a free $E$-module with basis
$x^{(a)}=x_1^{(a_1)}\cdot\ldots \cdot x_m^{(a_m)}$, $a\in \N^m$, $|a|=i$. The
$E$-linear differential on $C\mlpnt(v_1,\ldots,v_m;E)$ is
\begin{eqnarray*}
\partial_i : C_i(v_1,\ldots,v_m;E)&\longrightarrow& C_{i-1}(v_1,\ldots,v_m;E)\\
x^{(a)}&\mapsto& \sum_{a_j>0}v_jx_1^{(a_1)}\cdot\ldots\cdot x_j^{(a_j-1)}\cdot\ldots\cdot x_m^{(a_m)}.
\end{eqnarray*} 
One easily sees that $\partial\circ\partial=0$. So this is indeed a complex.

\begin{definition} 
Let $M \in \Mcc$.
The complex
$$
C\mlpnt(\vb;M)=C\mlpnt(\vb;E)\tensor_EM $$
is called the \emph{Cartan complex} of
$\vb$ with values in $M$. The corresponding homology modules
$$
H_i(\vb;M)=H_i(C\mlpnt(\vb;M)) $$
are called the \emph{Cartan homology} of $\vb$ with values in $M$.
\end{definition}

Cartan homology can be computed inductively as there is a long exact
sequence connecting the homologies of $v_1,\ldots,v_j$ and
$v_1,\ldots,v_j,v_{j+1}$ for $j=1,\ldots,m-1$.

To begin with there exists an exact sequence of complexes
$$0\longrightarrow C\mlpnt(v_1,\ldots,v_j;M)\stackrel{\iota}{\longrightarrow}C\mlpnt(v_1,\ldots,v_{j+1};M)\stackrel{\tau}{\longrightarrow}C\mlpnt(v_1,\ldots,v_{j+1};M)(-1)\longrightarrow 0,$$
where $(-1)$ indicates a shift in the homological degree, $\iota$ is the natural inclusion map and $\tau$ is given by 
$$\tau(g_0+g_1x_{j+1}+\ldots+g_kx^{(k)}_{j+1})=g_1+g_2x_{j+1}+\ldots+g_kx^{(k-1)}_{j+1},$$
where the $g_i$ belong to $C_{k-i}(v_1,\ldots,v_j;M)$.

This exact sequence induces a long exact sequence of homology modules.
\begin{proposition} \label{long cartan seq} \cite[Propositions 4.1]{AH}
Let
$M\in\Mcc$ and $\vb=v_1,\ldots,v_m \in E_1$. For all $1\leq j\leq m$
there exists a long exact sequence of graded $E$-modules
$$
\ldots\longrightarrow H_i(v_1,\ldots,v_j;M)\stackrel{\alpha_i}{\longrightarrow} H_i(v_1,\ldots,v_{j+1};M)\stackrel{\beta_i}{\longrightarrow} H_{i-1}(v_1,\ldots,v_{j+1};M)(-1)
$$
$$
\stackrel{\delta_{i-1}}{\longrightarrow} H_{i-1}(v_1,\ldots,v_j;M)\longrightarrow H_{i-1}(v_1,\ldots,v_{j+1};M)\longrightarrow \ldots
.$$
Here $\alpha_i$ is induced by the inclusion map $\iota$, $\beta_i$ by $\tau$ and $\delta_{i-1}$ is the conncting homomorphism, which acts as follows: if $z=g_0+g_1x_{j+1}+\ldots+g_{i-1}x^{(i-1)}_{j+1}$ is a cycle in $C_{i-1}(v_1,\ldots,v_{j+1};M)$, then $\delta_{i-1}(\overline{z})=\overline{g_0v_{j+1}}$.
\end{proposition}

Setting $\deg x_i=1$ induces a grading on the
complex and its homologies. 

\noindent The Cartan complex $C\mlpnt(v_1,\ldots,v_m;E)$
with values in $E$ is exact if the linear forms are $K$-linearly independent (see e.g., \cite[Remark3.4(3)]{AAH} or the proof of \cite[Theorem 2.2]{AHH}) and hence it is a minimal graded free
resolution of $H_0(v_1,\ldots,v_m;E)=E/(v_1,\ldots,v_m)$ over $E$. Thus it can be used to
compute $\Tor_i^E(E/(v_1,\ldots,v_m),\cdot)$.

\begin{proposition} \cite[Theorem 2.2]{AHH}\label{Tor and Cartan}
Let $M\in\Mcc$ and $\vb=v_1,\ldots,v_m\in E_1$ linearly independent over $K$. There are
isomorphisms of graded $E$-modules
$$
\Tor_i^E(E/(v_1,\ldots,v_m),M)\cong H_i(\vb;M) \qquad \text{ for all } i\geq 0.
$$
\end{proposition}

Recall that $H^j(M,v)$ for $j\in\Z$ denotes the $j$-th homology of the complex
\begin{equation}\label{simplicial homology}
(M,v): \qquad \ldots\longrightarrow M_{j-1}\stackrel{\cdot v}{\longrightarrow}M_j\stackrel{\cdot v}{\longrightarrow}M_{j+1}\longrightarrow \ldots\end{equation} 
and that $v$ is $M$-regular if and only if $H^j(M,v)=0$ for all $j\in\Z$.

We have now laid the required background in order to give the definition of the exterior annihilator numbers with respect to a certain sequence.

\begin{definition}
Let $v_1,\ldots,v_n$ be a basis of $E_1$ and let $M\in\Mcc$. The numbers
$$\alpha_{i,j}(v_1,\ldots,v_n;M)=\dim_KH^j(M/(v_1,\ldots,v_{i-1})M,v_i) $$
for $j\in\Z$ and $1\leq i\leq n$ are called the \emph{exterior annihilator numbers} of $M$ with respect to $v_1,\ldots,v_n$.
\end{definition}
So far, we have defined exterior annihilator numbers for an $E$-module which do depend on the chosen sequence. The following theorem justifies the definition of the so-called exterior generic annihilator numbers which are independent of the sequence.

\begin{theorem}\label{open set}
Let $J\subseteq E$ be a graded ideal. Then there exists a non-empty Zariski-open set $U\subseteq \GL_n(K)$ such that
$$\alpha_{i,j}(\gamma(e_1,\ldots,e_n);E/J)=\alpha_{i,j}(e_1,\ldots,e_n;E/\gin_{<_{\revlex}}(J))$$
for all $\gamma=(\gamma_{i,j})_{1\leq i,j\leq n}\in U$, where $\gamma(e_1,\ldots ,e_n)=(\gamma_{1,1}e_1+\ldots +\gamma_{n,1}e_n,\ldots,\gamma_{1,n}e_1+\ldots +\gamma_{n,n}e_n)$ and $<_{\revlex}$ denotes the reverse lexicographic order with respect to $e_1<\ldots<e_n$.
\end{theorem}

\begin{proof}
Let $U'=\{\varphi\in\GL_n(K)~:~\ini_{<_{\revlex}}(\varphi(J))=\gin_{<_{\revlex}}(J)\}$ be the non-empty Zariski-open set of linear transformations that can be used to compute the generic initial ideal of $J$. Set $U=\{\varphi^{-1}~:~\varphi\in U'\}$.
Let $\gamma=\varphi^{-1}\in U$ and set $v_i=\gamma(e_i)$ for $1\leq i\leq n$, i.e., $\varphi(v_i)=e_i$. As $\varphi$ is an automorphism, $E/\left(J+(v_1,\ldots,v_i)\right)$ and $E/\left(\varphi(J)+(e_1,\ldots,e_i)\right)$ have the same Hilbert function. \cite[Proposition 5.1]{AH} implies that $\ini_{<_{\revlex}}(\varphi(J)+(e_1,\ldots,e_i))=\gin_{<_{\revlex}}(J)+(e_1,\ldots,e_i)$ (observe that we use the reversed order on $[n]$). Therefore also $E/\left(J+(v_1,\ldots,v_i)\right)$ and $E/\left(\gin_{<_{\revlex}}(J)+(e_1,\ldots,e_i)\right)$ have the same Hilbert function.
The sequences
\begin{flushleft}
$0\longrightarrow H^j(E/\left(J+(v_1,\ldots,v_{i-1})\right),v_i)\longrightarrow\left(E/\left(J+(v_1,\ldots,v_i)\right)\right)_j$ 
\end{flushleft}
\begin{flushright}$\stackrel{\cdot v_i}{\longrightarrow}\left(E/\left(J+(v_1,\ldots,v_{i-1})\right)\right)_{j+1}
\longrightarrow\left(E/\left(J+(v_1,\ldots,v_{i})\right)\right)_{j+1}\longrightarrow 0$\end{flushright}
and
\begin{flushleft}$0\longrightarrow H^j(E/\left(\gin_{<_{\revlex}}(J)+(e_1,\ldots,e_{i-1})\right),e_i)\longrightarrow\left(E/\left(\gin_{<_{\revlex}}(J)+(e_1,\ldots,e_i)\right)\right)_j$
\end{flushleft}
\begin{flushright}$\stackrel{\cdot e_i}{\longrightarrow}\left(E/\left(\gin_{<_{\revlex}}(J)+(e_1,\ldots,e_{i-1})\right)\right)_{j+1}\longrightarrow\left(E/\left(\gin_{<_{\revlex}}(J)+(e_1,\ldots,e_{i})\right)\right)_{j+1}\longrightarrow 0$
\end{flushright}
are exact sequences of $K$-vector spaces (recall the definition of $H^j(-,-)$ from (\ref{simplicial homology})). The vector space dimensions of the three latter vector spaces in the two sequences coincide, hence it holds that 
\begin{eqnarray*}
\alpha_{i,j}(v_1,\ldots,v_n;E/J)&=&\dim_KH^j(E/\left(J+(v_1,\ldots,v_{i-1})\right),v_i)\\
&=& \dim_KH^j(E/\left(\gin_{<_{\revlex}}(J)+(e_1,\ldots,e_{i-1})\right),e_i)\\
&=& \alpha_{i,j}(e_1,\ldots,e_n;E/\gin_{<_{\revlex}}(J)).
\end{eqnarray*}
\end{proof}

As mentioned beforehand, now the following definition makes sense.

\begin{definition} \label{annihilators}
Let $J\subseteq E$ be a graded ideal. We set $$\alpha_{i,j}(E/J)=\alpha_{i,j}(e_1,\ldots,e_n;E/\gin_{<_{\revlex}}(J))$$
for $j\in\Z$ and $1\leq i\leq n$ and call these numbers the \emph{exterior generic annihilator numbers} of $E/J$.
\end{definition}

\begin{remark} \label{depth and annihilators}
 Let $J\subseteq E$ be a graded ideal. By definition of the $\alpha_{i,j}$ and the fact that $\gin_{<_{\revlex}}(\gin_{<_{\revlex}}(J))=\gin_{<_{\revlex}}(J)$, it holds that 
$$\alpha_{i,j}(E/J)=\alpha_{i,j}(E/\gin_{<_{\revlex}}(J)).$$
Set $\alpha_i(E/J)=\sum_{j\in\Z}\alpha_{i,j}(E/J)$ and $1\leq r\leq n$. Then $\alpha_{i}(E/J)=0$ for all $i\leq r$ if and only if $r\leq \depth_E (E/J)$.
This is an easy consequence of the fact that $e_1,\ldots,e_i$ is a regular sequence on $E/\gin_{<_{\revlex}}(J)$ if and only if $i\leq \depth_E (E/\gin_{<_{\revlex}}(J))=\depth_E(E/J)$ (see Lemma \ref{variables regular sequence}).
\end{remark}

It is well-known that being a regular sequence is a Zariski-open condition. One can prove it directly using that a sequence of linear forms is regular if and only if the first Cartan homology vanishes. But using the generic annihilator numbers provides a short proof.

\begin{proposition}\label{ext regular generic}
Let $J\subseteq E$ be a graded ideal and $\depth_EE/J=t$. Then there exists a non-empty Zariski-open set $U\subseteq GL_n(K)$ such that $\gamma_{1,1}e_1+\ldots +\gamma_{n,1}e_n,\ldots,\gamma_{1,t}e_1+\ldots +\gamma_{n,t}e_n$ is an $E/J$-regular sequence for all $\gamma=(\gamma_{i,j})_{1\leq i,j\leq n}\in U$. 
\end{proposition}
\begin{proof}
Let $U$ be the non-empty Zariski-open set as in Proposition \ref{open set}, i.e., such that the annihilator numbers with respect to sequences $v_1,\ldots,v_n$ induced by $U$ equal the generic annihilator numbers. Following Lemma \ref{variables regular sequence} $e_1,\ldots,e_t$ is a regular sequence on $E/\gin_{<_{\revlex}}(J)$ and therefore $$\alpha_{i,j}(v_1,\ldots,v_n;E/J)=\alpha_{i,j}(e_1,\ldots,e_n;E/\gin_{<_{\revlex}}(J))=0$$
for $i\leq t$. Thus $v_1,\ldots,v_t$ is regular on $E/J$.
\end{proof}

It is also well-known that the same statement is true over the polynomial ring. Nevertheless we were not able to give a reference for this fact. In \cite{Swartz} Swartz gives a proof for the special case of Stanley-Reisner rings of simplicial complexes. Therefore we include a short proof following ideas from Herzog using almost regular sequences in Section \ref{sect not minimal} where these are defined.

We can describe the numbers $\alpha_{i,j}$ as follows.

\begin{theorem}\label{description annihilators}
 Let $J\subseteq E$ be a graded ideal. Then
$$\alpha_{i,j}(E/J)=|\{\overline{e_F}\in E/\gin_{<_{\revlex}}(J)~:~\deg \overline{e_F}=j,\  \min F\geq i+1,\  \overline{e_F}\not=0,\  \overline{e_ie_F}=0\}|.$$ 
Here $\overline{e_F}$ denotes the projection of $e_F\in E$ on $E/\gin_{<_{\revlex}}(J)$.
\end{theorem}

\begin{proof}
Since $\alpha_{i,j}(E/J)=\alpha_{i,j}(E/\gin_{<_{\revlex}}(J))$ we may assume that $J=\gin_{<_{\revlex}}(J)$. Then  $\alpha_{i,j}(E/J)$ can be computed using the sequence $e_1,\ldots,e_n$, i.e., 
$$\alpha_{i,j}(E/J)=H^j(E/J+(e_1,\ldots,e_{i-1}),e_i)=\left(\frac{\left(J+(e_1,\ldots,e_{i-1})\right):e_i}{J+(e_1,\ldots,e_i)}\right)_j.$$
Thus
$$\alpha_{i,j}(E/J)=|\{e_F~:~\deg e_F=j, e_F\in(J+(e_1,\ldots,e_{i-1})):e_i, \ e_F\not\in J+(e_1,\ldots,e_i)\}|.$$
Let $e_F\in(J+(e_1,\ldots,e_{i-1})):e_i$ of degree $j$. Then $e_ie_F\in J+(e_1,\ldots,e_{i-1})$.
Since $J$ is a monomial ideal, either  $e_ie_F\in (e_1,\ldots,e_{i-1})$ or $e_ie_F\in J$. In the first case it follows that $e_F\in(e_1,\ldots,e_{i-1})$ so that we do not need to count it. In the second case, $e_F\not\in J+(e_1,\ldots,e_i)$ is equivalent to $e_F\not\in J$ and $e_F\not\in (e_1,\ldots,e_i)$ or equivalently $\overline{e_F}\not=0$ and $\min F\geq i+1$.
\end{proof}

In the special case of the exterior Stanley-Reisner ring of a simplicial complex $\Delta$, Theorem \ref{description annihilators} yields the following combinatorial description of the exterior generic annihilator numbers.

\begin{corollary}\label{description annihilators of complex}
Let $\Delta$ be a simplicial complex and let $\Delta^{e}$ be its exterior shifting. Then
$$\alpha_{i,j}(K\{\Delta\})=|\{F\in\Delta^e~:~ |F|=j,\  [i]\cap F=\emptyset,\  F\cup\{i\}\not\in\Delta^e\}|.$$
\end{corollary}

Using this description we are able to express the symmetric Betti numbers of the Stanley-Reisner ring of the exterior shifting of a simplicial complex as a linear combination of certain generic annihilator numbers.
One of the key ingredients for the proof of the formula aforementioned is the well-known Eliahou-Kervaire formula for the symmetric Betti numbers, for squarefree stable ideals. We recall this formula together with some definitions before giving our result.
For a mononomial ideal $I$ in the polynomial ring or in the exterior algebra we denote by $G(I)$ the unique minimal monomial generating set of $I$. Let further denote $G(I)_j$ the set of those monomials in $G(I)$ which are of degree $j$ for $j\in \Z$.
The Eliahou-Kervaire formula gives an explicit way of how one can compute the symmetric Betti numbers of a squarefree stable ideal.

\begin{proposition} \cite[Corollary 2.3]{AHH2} \label{Eliahou-Kervaire}
Let $I\subseteq S$ be a squarefree stable ideal. Then
$$\beta_{i,i+j}(S/I)=\sum_{u\in G(I)_{j+1}}\binom{n-\min(u)-j}{ i-1}.$$
\end{proposition}
For a simplicial complex $\Delta$, as $\Delta^e$ is shifted, the Stanley-Reisner ideal  $I_{\Delta^{e}}$ of the exterior shifting is a squarefree stable ideal. We need this property in the following. 

Our result expresses the annihilator numbers in terms of the minimal generators of $I_{\Delta^e}$ in the polynomial ring or of $J_{\Delta^e}$ in the exterior algebra.

\begin{proposition}\label{relation with polynomial ring}
Let $\Delta$ be a simplicial complex and $\Delta^e$ be its exterior shifting. Then
$$\alpha_{l,j}(E/J_{\Delta})=|\{u\in G(I_{\Delta^e})_{j+1}~:~\min(u)=l\}|=|\{u\in G(J_{\Delta^e})_{j+1}~:~\min(u)=l\}|.$$
In particular,
$$\beta^S_{i,i+j}(K[\Delta^e])=\sum_{l=1}^n \binom{n-l-j}{i-1}\alpha_{l,j}(K\{\Delta\}).$$
\end{proposition}

\begin{proof}
As shown in Corollary \ref{description annihilators of complex} the number $\alpha_{l,j}(E/J_{\Delta})$ counts the cardinality of the set 
$$\Ac=\{A\in\Delta^e~:~ |A|=j,\  [l]\cap A=\emptyset,\  A\cup\{l\}\not\in\Delta^e\}.$$
On the other hand the minimal generators of $I_{\Delta^e}$ or $J_{\Delta^e}$ are the monomials corresponding to minimal non-faces of $\Delta^e$, i.e. the elements of $\{u\in G(I_{\Delta^e})_{j+1}~:~\min(u)=l\}$ are the monomials $x_B$ such that $B$ lies in the set
$$\Bc=\{B\not\in\Delta^e~:~|B|=j+1,\ \min(B)=l,\  \partial(B)\subseteq\Delta^e\},$$
where $\partial(B)=\{F\subset B~:~ F\not=B\}$ denotes the boundary of $B$.

We show that there is a one-to-one correspondence between $\Ac$ and $\Bc$. 
Let $B\in \Bc$. Then $l\in B$ and $A=B\setminus\{l\}$ is an element in $\Ac$. 
Conversely if $A\in\Ac$ then $B=A\cup\{l\}\in\Bc$. The only non-trivial point here is to see that the boundary of $B$ is contained in $\Delta^e$. This holds since $\Delta^e$ is shifted and $\min(B)=l$.

The statement about the Betti numbers then follows from the Eliahou-Kervaire formula  for squarefree stable ideals (Proposition \ref{Eliahou-Kervaire}).
\end{proof}

The exterior generic annihilator numbers can also be used to compute the Betti numbers over the exterior algebra. This is analogous to a result over the polynomial ring, see \cite[Corollary 1.2]{CHH}. 
To this end we use the \emph{Cartan-Betti numbers} introduced by Nagel, R\"omer and Vinai in \cite{NRV}. 

\begin{definition}
Let $J\subseteq E$ be a graded ideal and let $v_1,\ldots,v_n$ be a basis of $E_1$. We set 
$$h_{i,j}(r)(v_1,\ldots,v_n;E/J)=\dim_KH_i(v_1,\ldots,v_r;E/J)_j,$$
where  $H_i(v_1,\ldots,v_n;E/J)$ denotes the $i$-th Cartan homology.
\end{definition}

Nagel, R\"omer and Vinai remarked that there exists a non-empty Zariski-open set $W$ such that the $h_{i,j}$ are constant on it. Therefore they define:

\begin{definition}
Let $J\subseteq E$ be a graded ideal and let $v_1,\ldots,v_n$ be a basis of $E_1$. We set 
$$h_{i,j}(r)(E/J)=h_{i,j}(r)(v_1,\ldots,v_n;E/J)$$
for $(v_1,\ldots,v_n)\in W$ as above and call these numbers the \emph{Cartan-Betti numbers} of $E/J$.
\end{definition}

For $r=n$, we obtain from Proposition \ref{Tor and Cartan} that the Cartan-Betti numbers of $E/J$ are the usual exterior graded Betti numbers of $E/J$, i.e., $h_{i,j}(n)(E/J)=\beta_{i,j}^E(E/J)$.

We formulate and prove the following result using the generic annihilator numbers. Plugging in the description of the generic annihilator numbers in terms of the minimal generators of $J_{\Delta^e}$ and taking into account that we use the reversed order on $[n]$, our result is the same as \cite[Theorem 2.4(i)]{NRV} which is a direct consequence of the construction of the Cartan homology for stable ideals in \cite[Proposition 3.1]{AHH}.

\begin{theorem}\label{CartanBetti}
Let $J\subseteq E$ be a graded ideal. Then 
$$h_{i,i+j}(r)\leq\sum_{k=1}^r \binom{r+i-k-1}{i-1}\alpha_{k,j}(E/J) \quad i\geq 1, j\geq 0$$
and equality holds for all $i\geq 1$ and $1\leq r\leq n$ if and only if $J$ is componentwise linear. 
\end{theorem}

\begin{proof}
Let $v_1,\ldots,v_n$ be a sequence of linear forms that can be used to compute the exterior generic annihilator numbers and the Cartan-Betti numbers of $E/J$ as well. (Such a sequence exists as both conditions are Zariski-open and the intersection of two non-empty Zariski-open set remains non-empty.) Set $\alpha_{i,j}=\alpha_{i,j}(E/J)$ and  
$$A_i= \KER \left(E/\left(J+(v_1,\ldots,v_i)\right)\stackrel{\cdot v_i}{\longrightarrow}E/\left(J+(v_1,\ldots,v_{i-1})\right)\right)$$
for $1\leq i\leq n$.
Then $A_i$ is a graded $E$-module and the $K$-vector space dimension of the $j$-th graded piece equals $\alpha_{i,j}(E/J)$. The above map occurs in the long exact sequence of Cartan homologies (see Proposition \ref{long cartan seq}) since the $0$-th Cartan homology is $$H_0(v_1,\ldots,v_r;E/J)=E/\left(J+(v_1,\ldots,v_r)\right).$$ 
Thus for $i=1$ and $r=1$ we obtain from the long exact Cartan homology sequence the following exact sequence
$$H_1(1)(-1)_{j+1}\rightarrow H_1(0)_{j+1}\rightarrow H_1(1)_{j+1}\rightarrow A_1(-1)_{j+1}\rightarrow 0.$$
Since $H_i(0)=0$ for $i\geq 1$ this yields
$$h_{1,j+1}(1)=\alpha_{1,j}.$$
For $r\geq 1$ we have the exact sequence
$$H_1(r+1)(-1)_{j+1}\rightarrow H_1(r)_{j+1}\rightarrow H_1(r+1)_{j+1}\rightarrow A_{r+1}(-1)_{j+1}\rightarrow 0.$$
From this sequence we conclude by induction hypothesis on $r$
\begin{eqnarray*}
 h_{1,j+1}(r+1)&\leq& \alpha_{r+1,j}+h_{1,j+1}(r)\\
&\leq& \alpha_{r+1,j}+\sum_{k=1}^r\binom{r-k}{0}\alpha_{k,j}\\
&=&\sum_{k=1}^{r+1}\alpha_{k,j}.\\
\end{eqnarray*}
Now let $i>1$. For $r=1$ there is the exact sequence
$$H_i(0)_{i+j}\rightarrow H_i(1)_{i+j}\rightarrow H_{i-1}(1)(-1)_{i+j}\rightarrow H_{i-1}(0)_{i+j}.$$
The outer spaces in the sequence are zero, hence 
$$h_{i,i+j}(1)=h_{i-1,i+j-1}(1)\leq\alpha_{1,j}$$
by induction hypothesis on $i$.

Let now $r\geq 1$. There is the exact sequence
$$H_i(r)_{i+j}\rightarrow H_i(r+1)_{i+j}\rightarrow H_{i-1}(r+1)(-1)_{i+j}\rightarrow H_{i-1}(r)_{i+j}.$$
We conclude by induction hypothesis on $r$ and $i$
\begin{eqnarray*}
 h_{i,i+j}(r+1)&\leq & h_{i,i+j}(r)+h_{i-1,i-1+j}(r+1)\\
&\leq& \sum_{k=1}^r\binom{r+i-k-1}{i-1}\alpha_{k,j}+\sum_{k=1}^{r+1}\binom{r+1+i-1-k-1}{i-2}\alpha_{k,j}\\
&=&\sum_{k=1}^r\left(\binom{r+i-k-1}{i-1}+\binom{r+i-k-1}{i-2}\right)\alpha_{k,j}+\binom{i-2}{i-2}\alpha_{r+1,j}\\
&=&\sum_{k=1}^{r+1}\binom{r+i-k}{i-1}\alpha_{k,j}.
\end{eqnarray*}

The inequalities in the proof are all equalities if and only if the long exact sequence is split exact.
In this case the sequence $v_1,\ldots,v_n$ is called a \emph{proper sequence} for $E/J$. In \cite[Theorem 2.10]{NRV} it is shown that this is the case if and only if $J$ is a componentwise linear ideal.
\end{proof}

\section{An unexpected behavior of the generic annihilator numbers}\label{sect not minimal}
A natural question to ask is whether the exterior generic annihilator numbers play a special role among the exterior annihilator numbers of $E/J$ with respect to a certain sequence. Herzog posed the question if they are the minimal ones among all the annihilator numbers. In the attempt of proving this conjecture it turned out to be wrong. In order to clarify this unexpected result we do not only give a counterexample of the conjecture but we also give a sketch of the original idea of the proof and explain how we came up with the example. This also gives a hint at how to construct further counterexamples. After some slight changes our example also serves as a counterexample of the corresponding conjecture for the symmetric generic annihilator numbers.
For the sake of completeness we first state the conjecture.

\begin{conjecture}
Let $J\subseteq E$ be a graded ideal. For any basis $v_1,\ldots,v_n$ of $E_1$ it holds that
$$\alpha_{i,j}(E/J)\leq \alpha_{i,j}(v_1,\ldots,v_n;E/J)$$
for $1\leq i\leq n$ and $j\geq 0$.
\end{conjecture}

Thus, our aim was to prove that the annihilator numbers are minimal on a non-empty Zariski-open set.
For $i=1$ this is known to be true. Just take the non-empty Zariski-open set such that the ranks of the matrices of the maps of the complex
$$(M,v): \qquad \ldots\longrightarrow M_{j-1}\stackrel{\cdot v}{\longrightarrow}M_j\stackrel{\cdot v}{\longrightarrow}M_{j+1}\longrightarrow \ldots$$
are maximal (note that $M_j=0$ for almost all $j$). To prove this for longer sequences we tried to show that the sets
$$U_{i,j}=\Big\{(v_1,\ldots,v_n)\subseteq E_1\mbox{ basis}~:~\substack{\alpha_{i,j}(v_1,\ldots,v_n;E/J)\leq \alpha_{i,j}(w_1,\ldots,w_n;E/J)\\ \mbox{for any basis } (w_1,\ldots,w_n)\subseteq E_1}\Big\}$$
were non-empty Zariski-open sets for $1\leq i\leq n$, $0\leq j\leq n$. The intersection of those sets would have been a non-empty Zariski-open set having the required property (Note that only finitely many sets are intersected.).
In order to compute a certain annihilator number $\alpha_{i,j}(v_1,\ldots,v_n;E/J)$ we used the exact sequence

\begin{flushleft}
$0\longrightarrow H^j(E/(J+(v_1,\ldots,v_{i-1})),v_i)\longrightarrow\left(E/\left(J+(v_1,\ldots,v_i)\right)\right)_j$ 
\end{flushleft}
\begin{flushright}$\stackrel{\cdot v_i}{\longrightarrow}\left(E/(J+(v_1,\ldots,v_{i-1}))\right)_{j+1}\longrightarrow\left(E/(J+(v_1,\ldots,v_{i}))\right)_{j+1}\longrightarrow 0$,
\end{flushright}

which yields
\begin{align*}
&\alpha_{i,j}(v_1,\ldots,v_n;E/J)=\dim_K\left(E/(J+(v_1,\ldots,v_i))\right)_j\\
&-(\dim_K\left(E/(J+(v_1,\ldots,v_{i-1}))\right)_{j+1}-\dim_K\left(E/(J+(v_1,\ldots,v_i))\right)_{j+1}).
\end{align*}
One can show that for a generic basis $\{v_1,\ldots,v_n\}\subseteq E_1$ each of the vector space dimensions
on the right-hand side of the above equality is minimized. However, being the left-hand side of the above equality an alternating sum of those three
dimensions there is no reason to expect it to be minimized by a generic basis.

\begin{example} \label{not minimal}
Let $1\leq i,j\leq n$ and let $J=(e_{l_1}\cdot \ldots \cdot e_{l_{j+1}}~:~i\leq l_1<l_2<\ldots <l_{j+1}\leq n)\subseteq E$ be a graded ideal. By construction, $J$ is strongly stable and it therefore holds that $\gin_{<_{\revlex}}(J)=J$. Then
$$\alpha_{i,j}(E/J)=\alpha_{i,j}(e_1,\ldots, e_n;E/\gin_{<_{\revlex}}(J))=\alpha_{i,j}(e_1,\ldots, e_n;E/J),$$
i.e., we can use the sequence $e_1,\ldots, e_n$ to compute the generic annihilator numbers of $E/J$.\\
From the exact sequence 
\begin{flushleft}
$0\longrightarrow H^j(E/(J+(e_1,\ldots,e_{i-1})),e_i)\longrightarrow\left(E/\left(J+(e_1,\ldots,e_i)\right)\right)_j$ 
\end{flushleft}
\begin{flushright}$\stackrel{\cdot e_i}{\longrightarrow}\left(E/(J+(e_1,\ldots,e_{i-1}))\right)_{j+1}\longrightarrow\left(E/(J+(e_1,\ldots,e_{i}))\right)_{j+1}\longrightarrow 0$
\end{flushright}
we deduce
\begin{align*}
\alpha_{i,j}=&\dim_K\left(E/(J+(e_1,\ldots,e_i))\right)_j\\
&-(\dim_K\left(E/(J+(e_1,\ldots,e_{i-1}))\right)_{j+1}-\dim_K\left(E/(J+(e_1,\ldots,e_i))\right)_{j+1}).
\end{align*}
In the following we consider the sequence ${\bf e}=e_1,\ldots,e_{i-2},e_i,e_{i-1},e_{i+1},\ldots,e_n$ and compute the exterior annihilator numbers of $E/J$ with respect to this sequence. As before, we have the exact sequence
\begin{flushleft}
$0\longrightarrow H^j(E/(J+(e_1,\ldots,e_{i-2},e_i)),e_{i-1})\longrightarrow\left(E/\left(J+(e_1,\ldots,e_i)\right)\right)_j$ 
\end{flushleft}
\begin{flushright}$\stackrel{\cdot e_{i-1}}{\longrightarrow}\left(E/(J+(e_1,\ldots,e_{i-2},e_i))\right)_{j+1}\longrightarrow\left(E/(J+(e_1,\ldots,e_{i}))\right)_{j+1}\longrightarrow 0$
\end{flushright}
which leads to
\begin{eqnarray*} 
\alpha_{i,j}({\bf e};E/J)&=&\dim_K\left(E/(J+(e_1,\ldots,e_i))\right)_j\\
&-&\!\!\!(\dim_K\left(E/(J+(e_1,\ldots,e_{i-2},e_i))\right)_{j+1}-\dim_K\left(E/(J+(e_1,\ldots,e_i))\right)_{j+1}).
\end{eqnarray*}
Our aim is to show that $\alpha_{i,j}(E/J)>\alpha_{i,j}({\bf e}; E/J)$.
We therefore need to show that $$\dim_K\left(E/(J+(e_1,\ldots,e_{i-2},e_i))\right)_{j+1}>\dim_K\left(E/(J+(e_1,\ldots,e_{i-2},e_{i-1}))\right)_{j+1}.$$
Let $m=e_{l_1}\cdot\ldots \cdot e_{l_{j+1}}\in E_{j+1}$ with $l_1<\ldots <l_{j+1}$. If $l_1\leq i-1$, it already holds that $m\in (e_1,\ldots,e_{i-1})_{j+1}\subseteq \left(J+(e_1,\ldots, e_{i-1})\right)_{j+1}$. If $l_1\geq i$, we have $i\leq l_1<\ldots<l_{j+1}$ and therefore $m\in J_{j+1}\subseteq (J+(e_1,\ldots, e_{i-1}))_{j+1}$. Thus, $m\in  (J+(e_1,\ldots, e_{i-1}))_{j+1}$ in either case and therefore $(J+(e_1,\ldots, e_{i-1}))_{j+1}=E_{j+1}$ and so $\dim_K\left(E/(J+(e_1,\ldots,e_{i-1}))\right)_{j+1}=0$.

Consider now $\widetilde{m}=e_{i-1}e_{i+1}\cdot\ldots\cdot e_{i+j}\in E_{j+1}$. By definition, it holds that $\widetilde{m}\notin J$ and $\widetilde{m}\notin (e_1,\ldots,e_{i-2},e_i)$. Since $J$ is a monomial ideal this implies $\widetilde{m}\notin (J+(e_1,\ldots,e_{i-2},e_i))_{j+1}$.
We therefore get $\dim_K\left(E/(J+(e_1,\ldots,e_{i-2},e_i))\right)_{j+1}>0$. This finally shows
$$\alpha_{i,j}(E/J)>\alpha_{i,j}({\bf e};E/J).$$

We also compute the annihilator numbers one step before w.r.t. this two sequences, i.e., $\alpha_{i-1,j}(E/J)$ and $\alpha_{i-1,j}({\bf e};E/J)$ to show that those numbers are related to each other the other way round, i.e., we have
\begin{equation}\label{second case}
\alpha_{i-1,j}(E/J)<\alpha_{i-1,j}({\bf e};E/J).
\end{equation}
This suggests that, to have a chance to become smaller than the generic numbers, the annihilator numbers with respect to ${\bf e}$ first have to become ``worse'', i.e., greater.
Similar to the $i$-th annihilator numbers of $E/J$ we can compute the $(i-1)$-st annihilator numbers using the exact sequence. We therefore get
\begin{align*}
\alpha_{i-1,j}(E/J)=&\dim_K\left(E/(J+(e_1,\ldots,e_{i-1}))\right)_j\\
&-(\dim_K\left(E/(J+(e_1,\ldots,e_{i-2}))\right)_{j+1}-\dim_K\left(E/(J+(e_1,\ldots,e_{i-1}))\right)_{j+1})
\end{align*}
for the $(i-1)$-st generic annihilator of $E/J$ in degree $j$. In the same way, we obtain
\begin{align*}
&\alpha_{i-1,j}({\bf e};E/J)=\dim_K\left(E/(J+(e_1,\ldots,e_{i-2},e_i))\right)_j\\
&-(\dim_K\left(E/(J+(e_1,\ldots,e_{i-2}))\right)_{j+1}-\dim_K\left(E/(J+(e_1,\ldots,e_{i-2},e_i))\right)_{j+1})
\end{align*}
for the $(i-1)$-st exterior annihilator number of $E/J$ in degree $j$ with respect to the sequence ${\bf e}$. Since $J$ is generated by monomials of degree strictly larger than $j$ it holds that
$$\dim_K\left(E/(J+(e_1,\ldots,e_{i-1}))\right)_j=\dim_K\left(E/(e_1,\ldots,e_{i-1})\right)_j$$
and 
$$\dim_K\left(E/(J+(e_1,\ldots,e_{i-2},e_i))\right)_j=\dim_K\left(E/(e_1,\ldots,e_{i-2},e_i)\right)_j.$$
Since
$$\dim_K\left(E/(J+(e_1,\ldots,e_{i-1}))\right)_j=\dim_K\left(E/(J+(e_1,\ldots,e_{i-2},e_i))\right)_j,$$
in order to show (\ref{second case}) we only need to prove that
$$\dim_K\left(E/(J+(e_1,\ldots,e_{i-1}))\right)_{j+1}<\dim_K\left(E/(J+(e_1,\ldots,e_{i-2},e_i))\right)_{j+1}.$$
This follows from 
\begin{equation}\label{containement}
\left(J+(e_1,\ldots,e_{i-1})\right)_{j+1}\supsetneq\left(J+(e_1,\ldots,e_{i-2},e_i)\right)_{j+1}.
\end{equation}
To show (\ref{containement}) let $m=e_{l_1}\cdot\ldots \cdot e_{l_{j+1}}\in \left(J+(e_1,\ldots,e_{i-2},e_i)\right)_{j+1}$ with $l_1<\ldots <l_{j+1}$. If $l_1\geq i$ it follows that $m\in J_{j+1}\subseteq \left(J+(e_1,\ldots,e_{i-1})\right)_{j+1}$. If $l_1\leq i-1$ it already holds that $e_{l_1}\in (e_1,\ldots,e_{i-1})$ and thus $m\in \left(J+(e_1,\ldots,e_{i-1})\right)_{j+1}$. Since $e_{i-1}e_{i+1}\cdot \ldots \cdot e_{i+j}\in \left(J+(e_1,\ldots,e_{i-1})\right)_{j+1}$ but $e_{i-1}e_{i+1}\cdot \ldots \cdot e_{i+j}\notin \left(J+(e_1,\ldots,e_{i-2},e_i)\right)_{j+1}$ we obtain (\ref{containement}).
\end{example}

After slight modifications to the above example (a kind of ``de-polarization'') we get a counterexample of the conjecture that the symmetric generic annihilator numbers are the minimal ones among all the annihilator numbers with respect to a sequence. Before treating this very example in more detail we recall the precise definition of the symmetric annihilator numbers. 
Let $v_1,\ldots,v_n$ be a $K$-basis of $S_1$ and let $M$ be a finitely generated graded $S$-module. We set $A_i(v_1,\ldots,v_n;M)=0:_{M/(v_1,\ldots,v_{i-1})M}v_i$ and call the numbers
\begin{align*}
\alpha_{i,j}(v_1,\ldots,v_n;M)=
\begin{cases}
\dim_K A_i(v_1,\ldots,v_n;M)_j, &\text{ if }1\leq i\leq n\\
\beta_{0,j}(M), &\text{ if } i=n+1\\
\end{cases}
\end{align*}
the {\itshape symmetric annihilator numbers} of $M$ with respect to the sequence $v_1,\ldots,v_n$.
Let now $v\in S_1$ be a linear form. If the multiplication map $M_{i-1}\rightarrow M_i$ is injective for all $i\gg 0$, i.e., if only finitely many graded components of $0:_M v$ are non-zero, $v$ is called an {\itshape almost regular element} on $M$.
A sequence $v_1,\ldots,v_r$ is an {\itshape almost regular sequence} on $M$ if $v_i$ is almost regular on $M/(v_1,\ldots,v_{i-1})M$ for $1\leq i\leq r$. Note that for an almost regular sequence only finitely many annihilator numbers are non-zero.

Herzog and Hibi proved in \cite{HH-book} that the set of almost regular sequences on $M$ is a non-empty Zariski-open set. They further show that $v_1,\ldots, v_t$ is an $M$-regular sequence if $v_1,\ldots,v_n$ is an almost regular sequence on $M$ and $\depth_S(M)=t$. Hence this gives a proof of the following well-known fact. 

\begin{proposition}
Let $M$ be a finitely generated graded $S$-module and let $\depth_S(M)=t$. There exists a Zariski-open set $U\subseteq GL_n(K)$ such that $\gamma_{1,1}x_1+\ldots +\gamma_{n,1}x_n,\ldots,\gamma_{1,t}x_1+\ldots +\gamma_{n,t}x_n$ is an $M$-regular sequence for all $\gamma=(\gamma_{i,j})_{1\leq i,j \leq n}\in U$. 
\end{proposition}

As for the exterior annihilator numbers Herzog and Hibi show that the symmetric annihilator numbers with respect to a sequence coincide when chosing the sequence from a certain non-empty Zariski-open set. The proof is very similar to the proof of Theorem \ref{open set}, the corresponding statement over the exterior algebra.
\begin{theorem} \cite{HH-book} 
Let $I\subseteq S$ be a graded ideal. Then there exists a non-empty Zariski-open set $U\subseteq GL_n(K)$ such that
$\gamma(x)=(\sum_{i=1}^n \gamma_{i,1}x_i,\ldots,\sum_{i=1}^n \gamma_{i,n} x_i)$ is almost regular for all $\gamma=(\gamma_{i,j})_{1\leq
i,j\leq n}\in U$. Moreover, 
$$\alpha_{i,j}(\gamma(x);S/I)=\alpha_{i,j}(x_1,\ldots,x_n;S/\gin_{<_{\revlex}}(I))$$
for all $\gamma\in U$.
\end{theorem}

This gives rise to the definition of the symmetric generic annihilator numbers. The numbers
$$\alpha_{i,j}(S/I)=\alpha_{i,j}(x_1,\ldots,x_n; S/\gin_{<_{\revlex}}(I))$$
are called the {\itshape symmetric generic annihilator numbers} of $S/I$.
As in the case of the exterior generic annihilator numbers one can wonder if the symmetric generic annihilator numbers are the minimal ones among the symmetric annihilator numbers with respect to a sequence.
As already mentioned, this is not the case. 

\begin{example}\label{not minimal 2}
Let $1\leq i\leq j\leq n$ and let $I=(x_{l_1}\cdot \ldots\cdot x_{l_{j+1}}~:~i\leq l_1\leq \ldots \leq l_{j+1})\subseteq S$ be a graded ideal. By construction $I$ is strongly stable and it therefore holds that $\gin_{<_{\revlex}}(I)=I$. Analogously to Example \ref{not minimal} we can use the sequence $x_1,\ldots, x_n$ to compute the symmetric generic annihilator numbers of $S/I$. From the exact sequence
\begin{flushleft}
$0\longrightarrow A_i(x_1,\ldots,x_n;S/I)_j\longrightarrow\left(S/\left(I+(x_1,\ldots,x_{i-1})\right)\right)_j$ 
\end{flushleft}
\begin{flushright}$\stackrel{\cdot x_i}{\longrightarrow}\left(S/(I+(x_1,\ldots,x_{i-1}))\right)_{j+1}\longrightarrow\left(S/(I+(x_1,\ldots,x_{i}))\right)_{j+1}\longrightarrow 0$
\end{flushright}
we deduce
\begin{align*}
\alpha_{i,j}(S/I)=&\dim_K\left(S/(I+(x_1,\ldots,x_{i-1}))\right)_j\\
&-(\dim_K\left(S/(I+(x_1,\ldots,x_{i-1}))\right)_{j+1}-\dim_K\left(S/(I+(x_1,\ldots,x_i))\right)_{j+1}).
\end{align*}
With the same exact sequence for the sequence ${\bf x}=x_1,\ldots,x_{i-2},x_i,x_{i-1},x_{i+1},\ldots,x_n$ we get
\begin{align*}
&\alpha_{i,j}({\bf x};S/I)=\dim_K\left(S/(I+(x_1,\ldots,x_{i-2},x_i))\right)_j\\
&-(\dim_K\left(S/(I+(x_1,\ldots,x_{i-2},x_i))\right)_{j+1}-\dim_K\left(S/(I+(x_1,\ldots,x_i))\right)_{j+1}).
\end{align*}
Since $I$ is generated in degree $j+1$ it holds that
\begin{align*}
\dim_K\left(S/(I+(x_1,\ldots,x_{i-1}))\right)_j&=\dim_K\left(S/(x_1,\ldots,x_{i-1})\right)_j\\
&=\dim_K\left( S/(x_1,\ldots,x_{i-2},x_i)\right)_j\\
&=\dim_K\left(S/(I+(x_1,\ldots,x_{i-2},x_i))\right)_j.
\end{align*}
One easily shows that $(I+(x_1,\ldots,x_{i-2},x_i))_{j+1}\subsetneq (I+(x_1,\ldots,x_{i-1}))_{j+1}$. This implies 
\begin{equation}\label{inequ: dimensions}
\dim_K\left(S/(I+(x_1,\ldots,x_{i-1}))\right)_{j+1}<\dim_K\left(S/(I+(x_1,\ldots,x_{i-2},x_i))\right)_{j+1}.
\end{equation}
As in the exterior case we thus obtain
$$\alpha_{i,j}(S/I)>\alpha_{i,j}({\bf x};S/I).$$
We now compute the $(i-1)$-st annihilator numbers in degree $j$ and see what happens in this case.
The same arguments as before show that 
\begin{align*}
\alpha_{i-1,j}(S/I)=&\dim_K\left(S/(I+(x_1,\ldots,x_{i-2}))\right)_j\\
&-(\dim_K\left(S/(I+(x_1,\ldots,x_{i-2}))\right)_{j+1}-\dim_K\left(S/(I+(x_1,\ldots,x_{i-1}))\right)_{j+1})
\end{align*}
and
\begin{align*}
&\alpha_{i-1,j}({\bf x};S/I)=\dim_K\left(S/(I+(x_1,\ldots,x_{i-2}))\right)_j\\
&-(\dim_K\left(S/(I+(x_1,\ldots,x_{i-2}))\right)_{j+1}-\dim_K\left(S/(I+(x_1,\ldots,x_{i-2},x_i))\right)_{j+1}).
\end{align*}
Using Equation (\ref{inequ: dimensions}) we thus obtain 
$$\alpha_{i-1,j}(S/I)<\alpha_{i-1,j}({\bf x};S/I).$$
\end{example}

Example \ref{not minimal} and \ref{not minimal 2} in particular show that changing the order of the elements of a sequence might change the annihilator numbers. However, when taking the sequence from a certain non-empty Zariski-open set the order of the elements does not matter.

\begin{theorem}\label{thm: order does not matter}
Let $J\subseteq E$ be a graded ideal. If $K$ has enough algebraically independent elements over its base field, then there exists a non-empty Zariski-open set $V\subseteq GL_n(K)$ such that
$$\alpha_{i,j}(E/J)=\alpha_{i,j}(\gamma(e_{\sigma(1)}),\ldots,\gamma(e_{\sigma(n)});E/J)$$
for all $\gamma \in V$ and all $\sigma\in S_n$, where $S_n$ denotes the symmetric group on $[n]$.
\end{theorem}
\begin{proof}
If $K$ has enough algebraically independent elements over its base field, then the subset $V\subseteq GL_n(K)$ consisting of matrices with algebraically independent entries is non-empty and Zariski-open. Furthermore, it is invariant under row permutations, column permutations and also under taking inverses. We show that $V$ is as required.

Consider $\gamma\in V$, $\sigma\in S_n$ and set $\gamma_i=\gamma(e_i)$.
The proof of Theorem \ref{open set} shows that it is enough to prove 
$$\dim_K \left( E/\left(\gin_{<_{\revlex}}(J)+(e_1,\ldots,e_i) \right)\right)_j=\dim_K \left( E/\left(J+(\gamma_{\sigma(1)},\ldots,\gamma_{\sigma(i)}) \right)\right)_j$$
for $i=1,\ldots,n$ and $j\in\Z$.

Since $\gin_{<_{\revlex}}(J)=\ini_{<_{\revlex}}(\gamma^{-1}(J))$ it holds that
$$\dim_K \left( E/\left(\gin_{<_{\revlex}}(J)+(e_1,\ldots,e_i) \right)\right)_j=\dim_K \left( E/\left(J+(\gamma_{1},\ldots,\gamma_{i}) \right)\right)_j.$$

The permutation lemma on page 288 in \cite{BK} (it is only formulated for monomial ideals, but the proof is valid also for arbitrary graded ideals) implies  that 

$$\dim_K \left( E/\left(J+(\gamma_{1},\ldots,\gamma_{i}) \right)\right)_j=\dim_K \left( E/\left(J+(\gamma_{\sigma(1)},\ldots,\gamma_{\sigma(i)}) \right)\right)_j.$$
Both equations together conclude the proof.
\end{proof}

It would be nice to know if one can drop the assumption on $K$ in the above Theorem \ref{thm: order does not matter}.

\section*{acknowledgement}
The main work of this article was made while the authors were participating at the summer school of PRAGMATIC 2008 in Catania, Sicily. The authors thank the organizers of PRAGMATIC 2008 for support and hospitality.
They are grateful to J\"urgen Herzog and Volkmar Welker for their suggestion to study this topic and for their helpful discussions. Further thanks go to Eran Nevo for pointing out the idea of Theorem \ref{characterization} to us and to the referee for very careful reading and for the proof of Theorem \ref{thm: order does not matter}.

\nocite{BS}
\bibliographystyle{amsalpha}
\bibliography{references_annihilators}

\providecommand{\bysame}{\leavevmode\hbox to3em{\hrulefill}\thinspace}
\providecommand{\MR}{\relax\ifhmode\unskip\space\fi MR }
\providecommand{\MRhref}[2]{%
  \href{http://www.ams.org/mathscinet-getitem?mr=#1}{#2}
}
\providecommand{\href}[2]{#2}
\begin{thebibliography}{AHH00}

\bibitem[AAH00]{AAH}
A.\ Aramova, L.~L.\ Avramov, and J.\ Herzog, \emph{Resolutions of monomial
  ideals and cohomology over exterior algebras}, Trans. Am. Math. Soc.
  \textbf{352} (2000), no.~2, 579--594.

\bibitem[AH00]{AH}
A.\ Aramova and J.\ Herzog, \emph{Almost regular sequences and {B}etti
  numbers}, Am. J. Math. \textbf{122} (2000), no.~4, 689--719.

\bibitem[AHH97]{AHH}
A.\ Aramova, J.\ Herzog, and T.\ Hibi, \emph{Gotzmann {T}heorems for {E}xterior
  {A}lgebras and {C}ombinatorics}, J. Algebra \textbf{191} (1997), 174--211.

\bibitem[AHH98]{AHH2}
\bysame, \emph{Squarefree lexsegment ideals}, Math. Z. \textbf{228} (1998),
  353--378.

\bibitem[AHH00]{AHH3}
\bysame, \emph{Shifting {O}perations and {G}raded {B}etti {N}umbers}, J.
  Alg.Comb. \textbf{12} (2000), 207--222.

\bibitem[BH98]{BH-book}
W.\ Bruns and J.\ Herzog, \emph{Cohen-{M}acaulay rings. {R}ev. ed.}, Cambridge
  Studies in Advanced Mathematics, vol.~39, Cambridge University Press, 1998.

\bibitem[BK88]{BK}
A.\ Björner and G.\ Kalai, \emph{{A}n extended {E}uler-{P}oincaré theorem},
  Acta Math. \textbf{161} (1988), no.~3-4.

\bibitem[BS87]{BS}
D.\ Bayer and M.\ Stillman, \emph{A criterion for detecting m-regularity},
  Invent. Math. \textbf{87} (1987), 1--11.

\bibitem[CHH04]{CHH}
A.\ Conca, J.\ Herzog, and T.\ Hibi, \emph{Rigid resolutions and big {B}etti
  numbers}, Comment. Math. Helv. \textbf{79} (2004), 826--839.

\bibitem[Eis04]{Ei}
D.\ Eisenbud, \emph{Commutative {A}lgebra with a {V}iew {T}oward {A}lgebraic
  {G}eometry}, Graduate Texts in Mathematics, vol. 150, Springer, 2004.

\bibitem[HH08]{HH-book}
J.\ Herzog and T.\ Hibi, \emph{Monomials}, Book manuscript, 2008.

\bibitem[HT99]{HT}
J.\ Herzog and N.\ Terai, \emph{Stable properties of algebraic shifting},
  Result in Math. \textbf{35} (1999), 260--265.

\bibitem[Kal01]{Kalai}
G.\ Kalai, \emph{Algebraic shifting}, Computational Commutative Algebra and
  Combinatorics. Adv. Stud. Pure Math. \textbf{33} (2001), 121--163.

\bibitem[KR09]{KR}
G.\ K\"ampf and T.\ R\"omer, \emph{Homological properties of {O}rlik-{S}olomon
  algebras}, Manuscr. Math. \textbf{129} (2009), 181--210.

\bibitem[KW08]{KW}
M.\ Kubitzke and V.\ Welker, \emph{The {M}ultiplicity {C}onjecture for
  {B}arycentric {S}ubdivisions}, Commun. Algebra \textbf{36} (2008), no.~11.

\bibitem[NRV08]{NRV}
U.~Nagel, T.~R\"omer, and N.~P. Vinai, \emph{Algebraic shifting and exterior
  and symmetric algebra methods}, Commun. Algebra \textbf{36} (2008), no.~1,
  208--231.

\bibitem[Röm01]{R}
T.\ Römer, \emph{Generalized {A}lexander {D}uality and {A}pplications}, Osaka
  J. Math. \textbf{38} (2001), 469--485.

\bibitem[Swa06]{Swartz}
E.\ Swartz, \emph{g-elements, finite buildings and higher {C}ohen-{M}acaulay
  connectivity}, J. Comb. Theory, Ser. A \textbf{113} (2006), no.~7,
  1305--1320.

\bibitem[Tru87]{Trung}
N.~V. Trung, \emph{Reduction exponent and degree bound for the defining
  equations of graded rings}, Proc. Amer. Math. Soc. \textbf{101} (1987),
  229--236.

\bibitem[Yan00]{Y}
K.\ Yanagawa, \emph{Alexander {D}uality for {S}tanley-{R}eisner {R}ings and
  {S}quarefree $\mathbb{N}^n$-{G}raded {M}odules}, J. Algebra \textbf{225}
  (2000), 630--645.

\end{thebibliography}

\end{document}